\documentclass[letterpaper, 10 pt, conference]{ieeeconf}  

\IEEEoverridecommandlockouts                              

\usepackage{graphicx,color} 
\usepackage{amsmath} 
\usepackage{amssymb}  

\title{\LARGE \bf
Optimization-based incentivization and control scheme for autonomous traffic
}

\author{Uro\v{s} Kalabi\'{c} \and Piyush Grover \and Shuchin Aeron
\thanks{This work was primarily supported by Mitsubishi Electric Research Laboratories. S.~Aeron performed some of this work while on an academic sabbatical with MERL in Spring 2019 and Summer 2019. S.~Aeron would like to acknowledge the support of NSF CAREER award CCF:1553075.}
\thanks{U.~Kalabi\'{c} is with Mitsubishi Electric Research Laboratories,
        Cambridge, MA 02139, USA
        {\tt\small kalabic@merl.com}}%
\thanks{P.~Grover is with the Department of Mechanical \& Materials Engineering, University of Nebraska, Lincoln, NE 68588, USA
        {\tt\small piyush.grover @unl.edu}}%
\thanks{S.~Aeron is with the Department of Electrical and Computer Engineering, Tufts University, Medford, MA 02155, USA
        {\tt\small shuchin @ece.tufts.edu}}%
}

\begin{document}

\maketitle
\thispagestyle{empty}
\pagestyle{empty}

\begin{abstract}
We consider the problem of incentivization and optimal control of autonomous vehicles for improving traffic congestion. In our scenario, autonomous vehicles must be incentivized in order to participate in traffic improvement. Using the theory and methods of optimal transport, we propose a constrained optimization framework over dynamics governed by partial differential equations, so that we can optimally select a portion of vehicles to be incentivized and controlled.

The goal of the optimization is to obtain a uniform distribution of vehicles over the spatial domain. To achieve this, we consider two types of penalties on vehicle density, one is the $L^2$ cost and the other is a multiscale-norm cost, commonly used in fluid-mixing problems. To solve this non-convex optimization problem, we introduce a novel algorithm, which iterates between solving a convex optimization problem and propagating the flow of uncontrolled vehicles according to the Lighthill-Whitham-Richards model. We perform numerical simulations, which suggest that the optimization of the $L^2$ cost is ineffective while optimization of the multiscale norm is effective. The results also suggest the use of a dedicated lane for this type of control in practice.


\end{abstract}

\section{Introduction}


Recent results show that the price of anarchy in vehicle traffic is essentially unbounded, \textit{i.e.}, traffic congestion can be arbitrarily worse as compared to the case when all vehicles on the road are routed by a central controller that aims to minimize congestion or average travel time \cite{zhang2016price,zhang2018price}.
The development of autonomous vehicles (AVs) is creating new possibilities for traffic management and control \cite{schwarting2018planning} and can be utilized for reducing the cost of anarchy. Keeping in mind that, in the autonomous future, vehicles will retain their primary function of transportation of individuals, we must consider that the desire of AV passengers may sometimes be antagonistic to the goals of traffic improvement. For example, given a choice of whether or not to let their vehicle participate in traffic improvement, passengers would likely consider the cost of their own time or comfort and, for this reason, would likely require some sort of incentive in order to participate.

In this work, we present an incentive and control scheme for the use of AVs to control traffic on a single road segment. We contrast this to most previous works that have considered incentivization for traffic management by improved routing on multiple road segments, \textit{e.g.}, \cite{mehr2019will,mehr2019pricing,tanaka2019linearly}; we instead focus on dynamic control of vehicles. Specifically, we model traffic as a flow and assume that a certain part of this flow can be incentivized to switch and cooperate in order to control traffic. This turns the problem of controlling traffic into both an incentivization and control problem. AV passengers must first be incentivized to cooperate, and then their vehicles must be controlled. We note that a related mixed traffic model using mean-field ideas was recently introduced in \cite{huang2019stabilizing}, however the focus in that work is on stability and not on incentivization.




The problem of computing appropriate incentives falls under the umbrella of mechanism design \cite{borgers2015introduction}. We are interested in developing an incentivization framework in the context of mixed-traffic control, where both cooperative and non-cooperative vehicles are present. As a first step towards this goal, in this work we consider the problem of selecting a subset of vehicles, periodically updated, that can be potentially incentivized, and hence optimally controlled to drive the overall vehicle density to be close to uniform. As part of the solution, the corresponding optimal controls for this subset are also computed. This problem is generally motivated by recent results \cite{wu2018stabilizing,wu2017emergent} showing that intelligent control of a small fraction of vehicles can significantly reduce congestion and attenuate traffic waves.

We propose a periodically updated, receding-horizon control scheme. During each period, we determine the distribution of incentive, \textit{i.e.}, which subset of the vehicles to control, as well as the control inputs to these vehicles over the fixed time-period. The receding horizon approach is used to address undesirable edge effects that typically appear at the end of the control time-horizon. The control history is determined according to the assumption that controlled AV flow can be directly controlled and that uncontrolled flow follows the conventional Lighthill-Whitham-Richards (LWR) model. The resulting problem is a variation of the mean-field optimal control problem \cite{Bensoussan2013} with congestion, where only a subset of the agents are controlled.

We consider two types of cost on vehicle density. The first is an $L^2$ cost, which penalizes the distance to the uniform density, and the second is an $\dot H^{-1}$ cost, which is a type of multiscale norm and used in fluid mixing to ensure uniformity.
To show how one can implement our scheme,
we perform a numerical simulation on a single, circular road scenario and use the optimization software \texttt{cvx} to solve the optimization problem. Our results show that the penalization of $L^2$ cost does not effectively achieve convergence to a uniform distribution while penalization of $\dot H^{-1}$ cost does. For this reason, we suggest that approaches from fluid mixing can be helpful in solving traffic optimization. Notably, we find that the scheme requires flow rates that are too high to achieve in a single lane; we therefore suggest that the practical implementation of our approach may require use of a dedicated lane.

The rest of the paper is organized as follows. In Section \ref{sec:prob}, we present a detailed problem formulation and introduce the optimization framework. In Section \ref{sec:algo}, we present a computationally feasible algorithm for solving the proposed optimization problem. In Section \ref{sec:sims}, we present detailed simulation results and provide a discussion on how one may implement the resulting scheme in practice. Section \ref{sec:conclude} is the conclusion.


\section{Problem Formulation}
\label{sec:prob}

We consider a two-population traffic scenario, consisting of heteronomous (non-cooperative) vehicles (HVs) and autonomous (cooperative) vehicles (AVs). The densities of HVs and AVs are given by $\rho_1(x,t), \rho_2(x,t) \in \mathbb R$, respectively, where $(x,t) \in \Omega \times [0,T]$ is the position $x$ of a vehicle at time $t$, and $\Omega$ and $[0,T]$ are the spatial and time domains, respectively. The vehicular dynamics of HVs are governed by the Lighthill-Whitham-Richards (LWR) model,\footnote{While in this paper we only report results for the LWR model, we note that it is straightforward to extend our algorithm to other models of uncontrolled traffic flow such as the Payne-Whitham and Aw-Rascle-Zhang (ARZ) models.} with the dynamics given by,
\begin{equation}\label{equ:HV_dyn}
\partial_t \rho_1 + \nabla \cdot \left(\rho_1 v_1(\rho_1,\rho_2)\right) = 0,
\end{equation}
where,
\begin{equation}\label{equ:LWR}
v_1(\rho_1,\rho_2) = u_0 \left( 1 - \frac{\rho_1 + \rho_2}{{\rho}^*}\right),
\end{equation}
is the HV velocity field and the parameters $u_0$ and $\rho^*$ are the free-flow speed and the maximum density, respectively. The density evolution of AVs is governed by the continuity equation,
\begin{equation}\label{equ:AV_dyn}
\partial_t\rho_2 + \nabla \cdot m_2(\rho_2,v_2) = 0,
\end{equation}
where $m_2(x,t) \in \mathbb R$ is the controlled flux, given by,
\begin{equation*}
m_2(\rho_2,v_2) = \rho_2v_2,
\end{equation*}
where $v_2(x,t) \in \mathbb R$ is the controlled, AV velocity field. Note that the HV dynamics are affected by the AV dynamics but not \textit{vice versa}.

Apart from being able to control the AV velocity field, we are able to assign the initial densities of AVs $\rho_2(x,0)$ over $\Omega$. We do this by assuming that a certain distribution of vehicles can be either an HV or AV, but must be incentivized to turn on AV capability. We assume that the total amount available for incentivization is fixed and an offer can be made only once over a time period $[0,T]$. As a consequence of optimality, there is no reason to not make an offer after the initial time $t=0$ since, if $v_2(x,t) = v_1(x,t)$ is optimal over $t \in [0,a]$ but not over $t \in [a,T]$, it will be equivalent to make an offer at any time $t \leq a$.
In this work, we assume that each vehicle operator's cost preference $\beta(x)$ is equal, \textit{i.e.}, $\beta(x)$ is constant over $x$. 

We assume the initial density distribution $\rho_0 = \rho_1(\cdot,0)+\rho_2(\cdot,0)$ is given and its average is $\bar \rho$, \textit{i.e.}, $\int_{\Omega} \rho_0(x)\text{d}x = \bar \rho$.
We also assume that the initial density distribution of potential AVs $\hat\rho_0 \leq \rho_0$, with average $\hat{\bar\rho}$, is given and
therefore the initial distribution of AVs $\rho_{2,0} = \rho_2(\cdot,0)$ must satisfy the incentivization condition,
\begin{equation}\label{equ:rho2_constr}
\int_{\Omega} \rho_2(x,0)\text{d}x \leq \int_{\Omega} \beta(x)\hat\rho_0(x)\text{d}x = \beta \hat{\bar\rho}.
\end{equation}

\subsection{Optimization}

Given an initial density distribution $\rho_0(x)$, we wish to determine the initial AV distribution $\rho_{2,0}: \Omega \to \mathbb R$ satisfying \eqref{equ:rho2_constr} and a controlled velocity profile $v_2: \Omega \times [0,T] \to \mathbb R$ to transport the density from $\rho_0$ to a density that is close to uniform over $\Omega$.

We consider two approaches in order to achieve this end. 
Both minimize a cost function of the form,
\begin{equation}
\iint  \rho_2(x,t)v_2(x,t)^2\text{d}x\text{d}t + cV(\rho_1,\rho_2), 
\end{equation}
so that the first term in the weighted sum is a cost on control and the second term is a cost on the state, and $c > 0$ is a scalar weight.

In the first case, the cost function is chosen with the aim of minimizing the $L^2$ cost of the density over $\Omega \times [0,T]$,
\begin{equation*}
V(\rho_1,\rho_2) = \iint (\rho_1(x,t)+\rho_2(x,t))^2 \text{d}x\text{d}t.
\end{equation*}
which is equivalent to minimizing the distance to the uniform density at all times. By performing the Benamou-Brenier change of variables \cite{benamou2000computational} $(\rho_2,v_2) \mapsto (\rho_2,m_2)$, the optimization problem can be written as,
\begin{equation}\label{equ:L2_cost}
\min_{\rho_{2,0}, m_2} \iint \frac{m_2^2}{\rho_2} \text{d}x\text{d}t + c\iint (\rho_1+\rho_2)^2 \text{d}x\text{d}t,
\end{equation}
subject to the constraints,
\begin{subequations}\label{equ:opt_constr}
\begin{align}
\rho_1(x,t),~\rho_2(x,t),~m_2(x,t) &\geq 0, \\
\rho_1(x,t) + \rho_2(x,t) &\leq \rho^*, \\
\int_{\Omega} \rho_{2,0}(x)\text{d}x &\leq \beta\hat{\bar\rho}, \\
\rho_1(x,0)+\rho_{2,0}(x) &= \rho_0(x),
\end{align}
\end{subequations}
for all $(x,t) \in \Omega \times [0,T]$, and dynamics \eqref{equ:HV_dyn}-\eqref{equ:AV_dyn}.

In the second case, the cost function is chosen with the aim of minimizing the square of the $\dot H^{-1}$ seminorm \cite{multiscale_norms} of the density over $\Omega \times [0,T]$. This cost function is used in fluid mixing problems \cite{multiscale_norms} to optimize mixing. Setting a penalty on the deviation from the mean, the cost is given by,
\begin{multline*}
V(\rho_1,\rho_2) \\ = \frac{1}{\bar\rho T}\iint \left|(-\Delta)^{-\frac{1}{2}}(\rho_1(x,t)+\rho_2(x,t))\right|^2\text{d}x\text{d}t-1.
\end{multline*}
The cost can be reformulated in terms of Fourier coefficients to obtain the full optimization problem,
\begin{equation}
\min_{\rho_{2,0}, m_2} \iint \frac{m_2^2}{\rho_2} \text{d}x\text{d}t + \frac{c}{\bar\rho}\int\sum_{k=2}^\infty \frac{(\hat\rho_{1,k}+\hat\rho_{2,k})^2}{k^2} \text{d}t,\label{equ:H1_cost}
\end{equation}
subject to the constraints \eqref{equ:opt_constr}, where $\hat\rho_{1,k}$ and $\hat\rho_{2,k}$ are the $k$-th Fourier coefficients of $\rho_1(\cdot,t)$ and $\rho_2(\cdot,t)$, respectively. Hence, this cost weighs the spatial low-frequency components of fluctuations in the density higher than the high-frequency components. It has been shown to have good numerical properties in practical optimization problems \cite{foures2014optimal}.

The optimization problem can be considered as a variation of a mean-field optimal control problem \cite{Bensoussan2013} with congestion, where only a subset of the agents are directly controlled. Since the problem is non-convex due to the LWR dynamics in \eqref{equ:LWR}, the conventional approaches to solving convex mean-field control problems do not apply, and we must design a new algorithm.

\section{Algorithm}
\label{sec:algo}

In both cases above, the nonlinear dynamics \eqref{equ:HV_dyn}-\eqref{equ:AV_dyn} result in non-convex problems. This is because the flux expression obtained by multiplying \eqref{equ:LWR} by the density $\rho_1$ cannot be converted to a linear PDE. We can see this when we discretize the dynamics using a conservative scheme, like that of Godunov, where the flux term becomes dependent on logical conditions.

The problem can be somewhat alleviated by using a discretization of the PDEs that is dissipative, such as that of Lax-Friedrichs. In this way, the map from the state $(\rho_1(\cdot,t_k),\rho_2(\cdot,t_k))$ at one time instant $t_{k}$ to the next time instant $t_{k+1}$ becomes linear. The drawback, however, is that the use of a dissipative discretization does not lead to a physically correct solution \cite{leveque_book}.

We propose to use a two step iterated approach  in solving the optimization problem.
\begin{itemize}
    \item {Step 1:}
Fix the density $\rho_1$ of the HVs (in space and time), and use the Lax-Friedrichs discretization to obtain an optimized control $v_2$ and AV density $\rho_2$.

\item{Step 2:} Propagate the HV dynamics using a first-order Godunov scheme to determine a new value for HV density, using the AV density obtained in step 1. 
\end{itemize}
In this way, by keeping $\rho_1$ static in step 1, the corresponding optimization problem becomes convex and the dissipative components of Lax-Friedrichs can be compensated through vehicle actuation. 

At the beginning of the algorithm, we set initial conditions $\rho_1^{(0)}$ and $\rho_2^{(0)}$ that satisfy the constraints \eqref{equ:opt_constr}.
Then, during each iteration $i \geq 1$, we perform the two steps. 

In the first step, we solve the optimization problem,
\begin{equation}\label{equ:opt2_cost}
\min_{\rho_{2,0}, m_2} \iint \frac{m_2(x,t)^2}{\rho_2(x,t)} \text{d}x\text{d}t + cV\left(\rho_1^{(i-1)},\rho_2\right),
\end{equation}
subject to the constraints,
\begin{subequations}\label{equ:opt2_constr}
\begin{align}
\rho_2(x,t),~m_2(x,t) &\geq 0, \\
\rho_2(x,t) &\leq \bar\rho - \rho_1^{(i-1)}(x,t), \\
\int_{\Omega} \rho_{2,0}(x)\text{d}x &\leq \beta\hat{\bar\rho}, \\
\rho_{2,0}(x) &= \rho_0(x)-\rho_1^{(i-1)}(x,0),
\end{align}
\end{subequations}
for all $(x,t) \in \Omega \times [0,T]$, and dynamics \eqref{equ:AV_dyn}. Upon discretization, the dynamics \eqref{equ:AV_dyn} are linear, of the form,
\begin{equation}
\rho_2(x_k,t_{k+1}) =
D_\text{LF}\begin{bmatrix}
\rho_2(\cdot,t_{k}) \\ m_2(\cdot,t_k)\end{bmatrix},
\end{equation}
where $D_\text{LF}$ is the fixed, spatial discretization matrix obtained using the Lax-Friedrichs method and $t_k$ is the $k$-th time step.

In the second step, we use a first-order Godunov method to solve the dynamic equation,
\begin{equation}\label{equ:HV_dyn2}
\partial_t \rho_1 + u_0\nabla \cdot \rho_1\left(1 - \frac{\rho_1 + \rho_2^{(i)}}{\rho^*}\right) = 0,
\end{equation}
with initial condition $\rho_1(x,0) = \rho_0(x)-\rho_{2,0}^{(i)}(x)$,
alternating between steps until either a stopping criterion or the maximum number of iterations have been reached.

\subsubsection{Convergence}
The problem \eqref{equ:opt2_cost}-\eqref{equ:opt2_constr} is convex and the solution to the dynamics \eqref{equ:HV_dyn2} is unique due to continuity \cite{modelingbook16}. These are desirable properties, but do not guarantee convergence. The constraints prevent the algorithm from blowing up but do not prevent limit cycles from occuring. For this reason, we terminate the algorithm after a number of prescribed iterations.

\subsection{Receding-Horizon Approach}
To ensure that the control scheme operates as seamlessly as possible, and because optimizing over a finite time-horizon can often result in undesirable effects as the time reaches the end of the optimization, we consider the use of a receding-horizon approach.

Specifically, we solve the problem \eqref{equ:opt2_cost}-\eqref{equ:HV_dyn2} over the time-interval $[0,T]$ but only implement the control over $[0,T/N]$ where $N$ is a design parameter. We then solve the same problem over the time-interval $[T/N,T+T/N]$ with initial conditions set to the values obtained at time $T/N$ and implement the control over $[T/N,2T/N]$. We repeat this procedure for times $2T/N,3T/N,\dots$ until we obtain a solution over the entire time-interval $[0,T]$.

\begin{table}[t]
\caption{Simulation Parameters}
\label{tab:params}
\begin{center}
\begin{tabular}{|c|c|c|}
\hline
Parameter & Value & Units \\
\hline
$|\Omega|$ & $2$ & mi \\
$T$ & $8$ & min \\
$u_0$ & $1$ & mi/min \\
$\rho^*$ & $10$ & $50$ veh/mi \\
$c$ & $0.1$ & \\
\hline
\end{tabular}
\end{center}
\end{table}

\section{Numerical Simulations}
\label{sec:sims}

\begin{figure}[t]
\begin{center}
\includegraphics[width=0.45\textwidth]{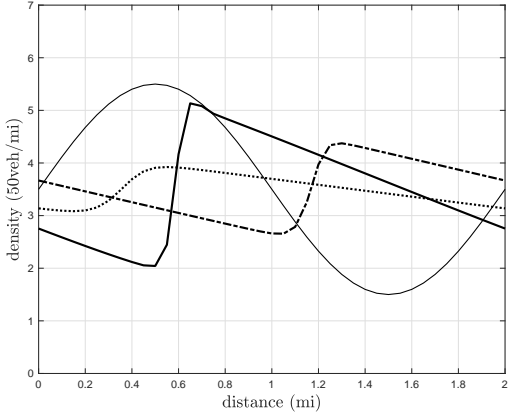}    
\caption{Initial distribution of vehicles $\rho_0$ (solid) and average (dotted)} 
\label{fig:rho0}
\end{center}
\end{figure}

We perform numerical simulations to test and exhibit the properties of the algorithm just presented. In all simulations, the simulation parameters are set to the values given in Table \ref{tab:params} and the optimization problems are solved using the software package \texttt{cvx}. We begin with the initial distribution $\rho_0$ over $\Omega = S^1$ given in Fig.~\ref{fig:rho0} as,
\begin{equation*}
\rho_0(x) = 3.5 + 2\sin(\pi x).
\end{equation*}
For comparison, the figure also shows the evolution of vehicle density without control, \textit{i.e.}, $v_2 = v_1(\rho_1,\rho_2)$. We begin by considering the optimization of the $L^2$ cost on density and then the $\dot H^{-1}$ cost.

\begin{figure}[t]
\begin{center}
\includegraphics[width=0.45\textwidth]{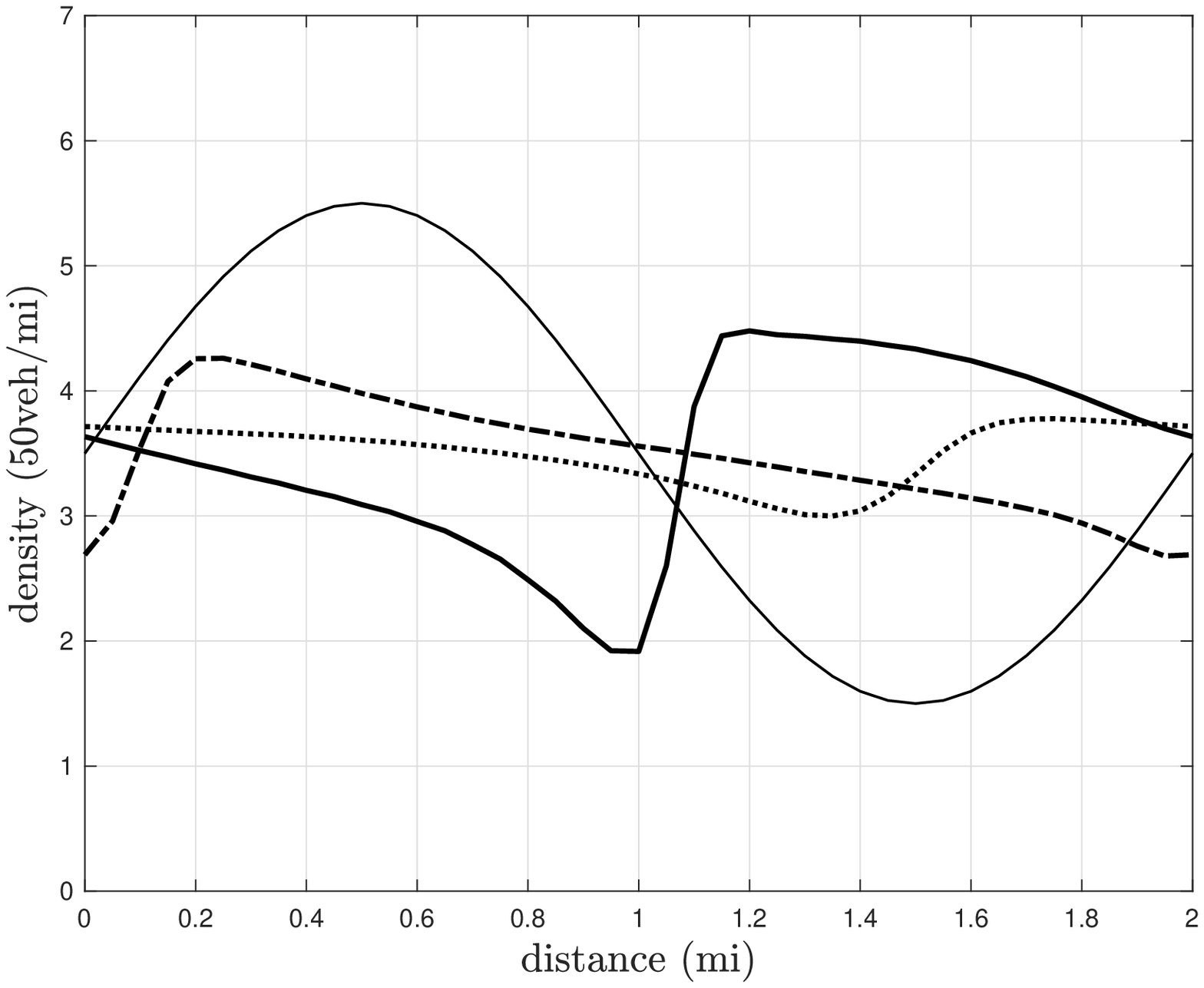}    
\caption{Total vehicle density at $t = 0$ (solid, light), $t = T/4$ (solid), $t = T/2$ (dot-dashed), and $t = T$ (dotted) corresponding to optimization of $L^2$ cost with $N=1$} 
\label{fig:rho_hist_L2_N1}
\end{center}
\end{figure}

\begin{figure}[t]
\begin{center}
\includegraphics[width=0.45\textwidth]{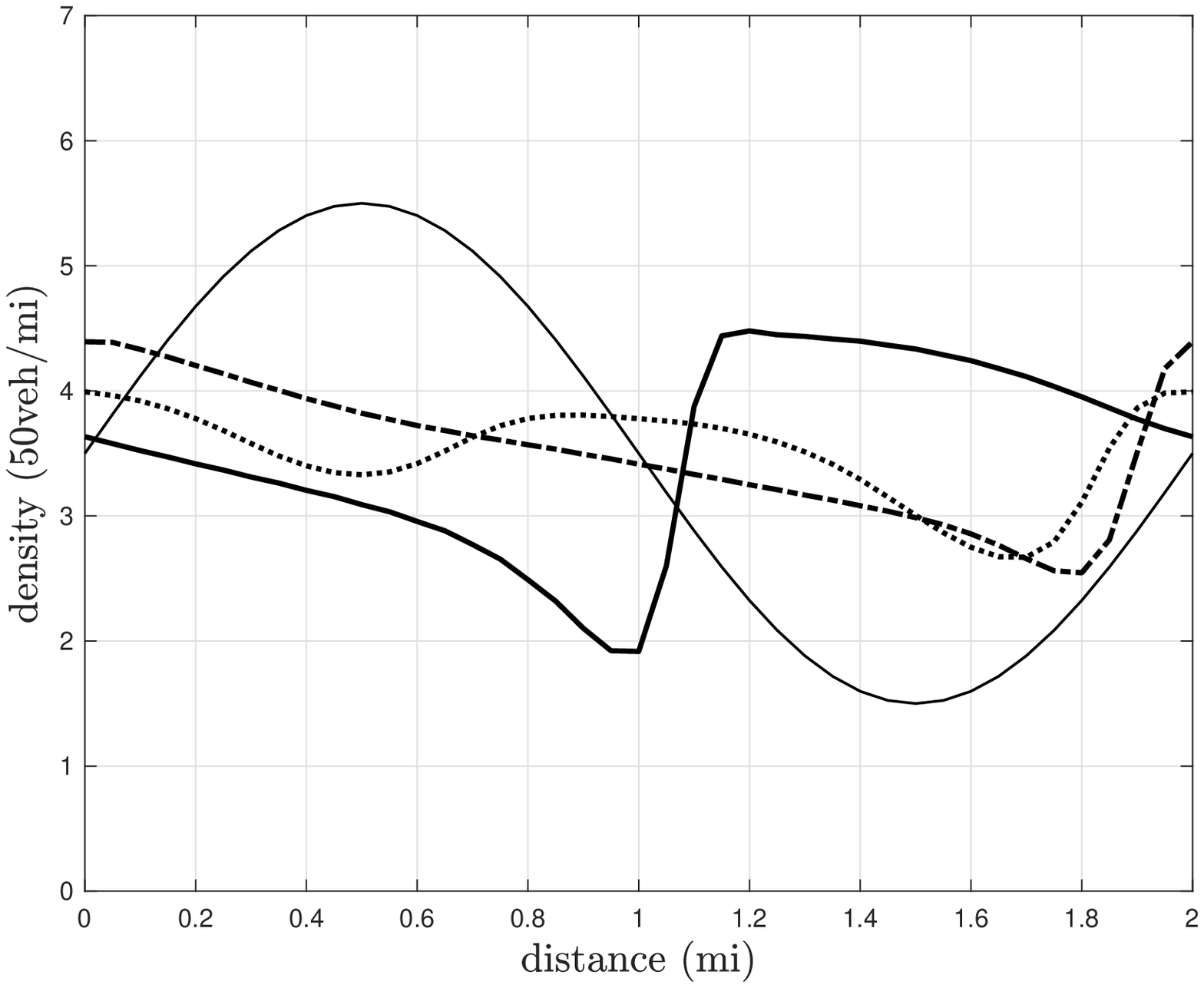}    
\caption{Total vehicle density at $t = 0$ (solid, light), $t = T/4$ (solid), $t = T/2$ (dot-dashed), and $t = T$ (dotted) corresponding to optimization of $L^2$ cost with $N=2$} 
\label{fig:rho_hist_L2_N2}
\end{center}
\end{figure}

\begin{figure}[t]
\begin{center}
\includegraphics[width=0.45\textwidth]{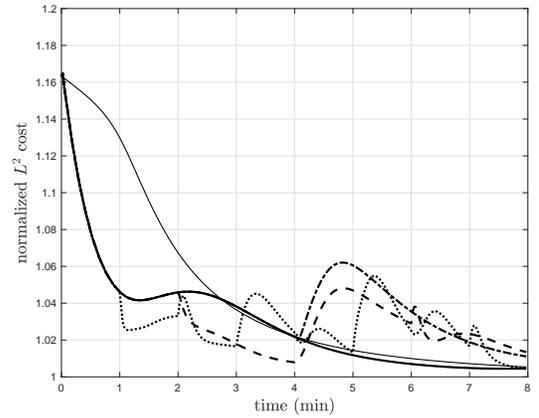}    
\caption{$L^2$ cost density component corresponding to uncontrolled (solid, light), $N=1$ (solid), $N=2$ (dot-dashed), $N=4$ (dashed), and $N=8$ (dotted)}
\label{fig:L2_cost}
\end{center}
\end{figure}


\begin{figure}[t]
\begin{center}
\includegraphics[width=0.45\textwidth]{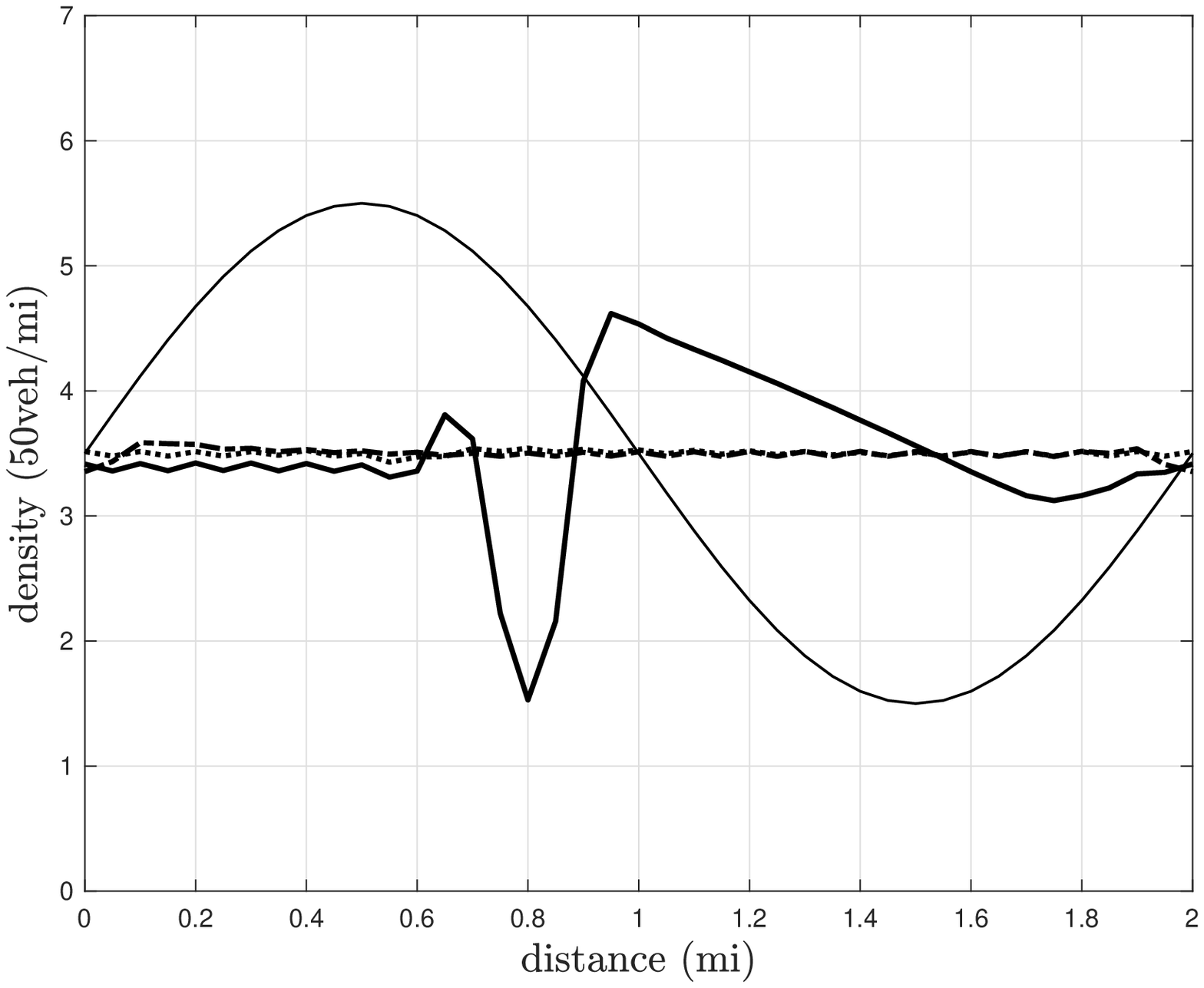}    
\caption{Total vehicle density at $t = 0$ (solid, light), $t = T/4$ (solid), $t = T/2$ (dot-dashed), and $t = T$ (dotted) corresponding to optimization of $\dot H^{-1}$ cost with $N=2$} 
\label{fig:rho_hist_H1_N2}
\end{center}
\end{figure}

\begin{figure}[t]
\begin{center}
\includegraphics[width=0.45\textwidth]{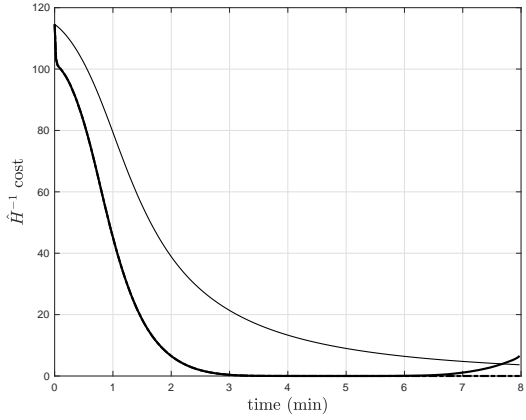}    
\caption{$\dot H^{-1}$ cost density component corresponding to uncontrolled (solid, light), $N=1$ (solid), and $N=2$ (dot-dashed)}
\label{fig:H1_cost}
\end{center}
\end{figure}


\begin{figure}[t]
\begin{center}
\includegraphics[width=0.45\textwidth]{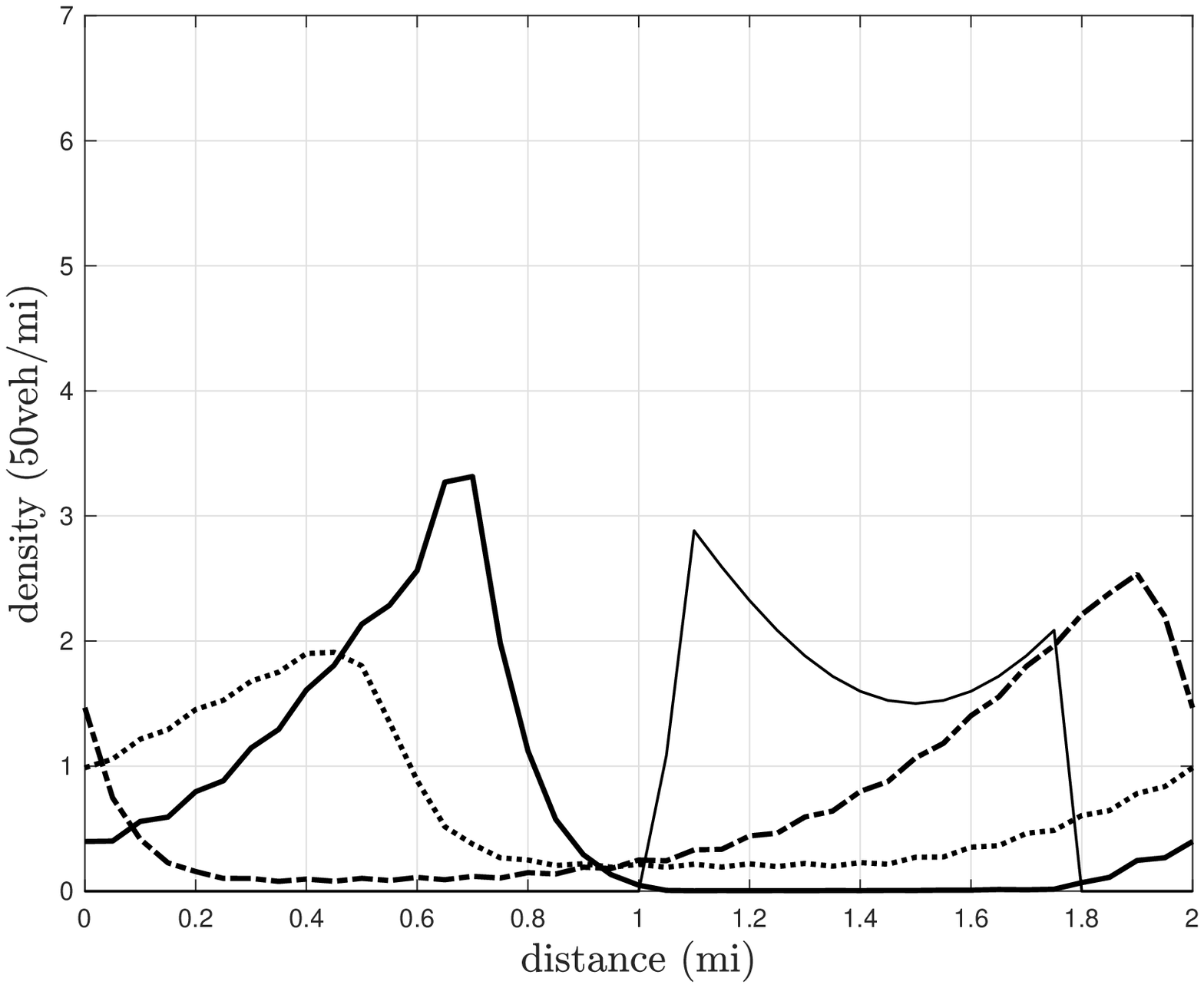}    
\caption{AV density at $t = 0$ (solid, light), $t = T/4$ (solid), $t = T/2$ (dot-dashed), and $t = T$ (dotted) corresponding to optimization of $\dot H^{-1}$ cost with $N=2$} 
\label{fig:rho2_hist_H1_N2}
\end{center}
\end{figure}

\begin{figure}[t]
\begin{center}
\includegraphics[width=0.45\textwidth]{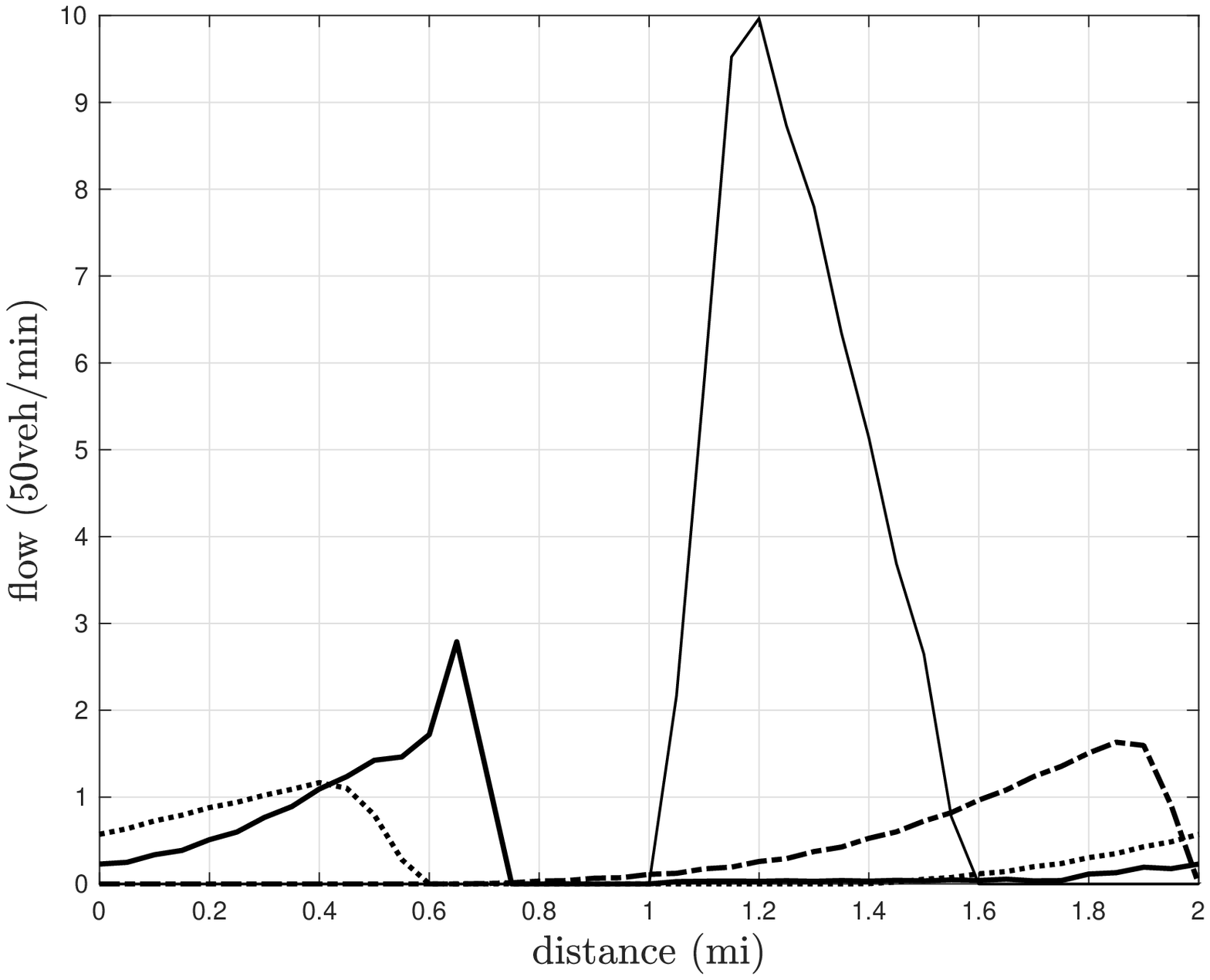}    
\caption{AV flow at $t = 0$ (solid, light), $t = T/4$ (solid), $t = T/2$ (dot-dashed), and $t = T$ (dotted) corresponding to optimization of $\dot H^{-1}$ cost with $N=2$} 
\label{fig:H1_flow}
\end{center}
\end{figure}

\begin{figure}[t]
\begin{center}
\includegraphics[width=0.45\textwidth]{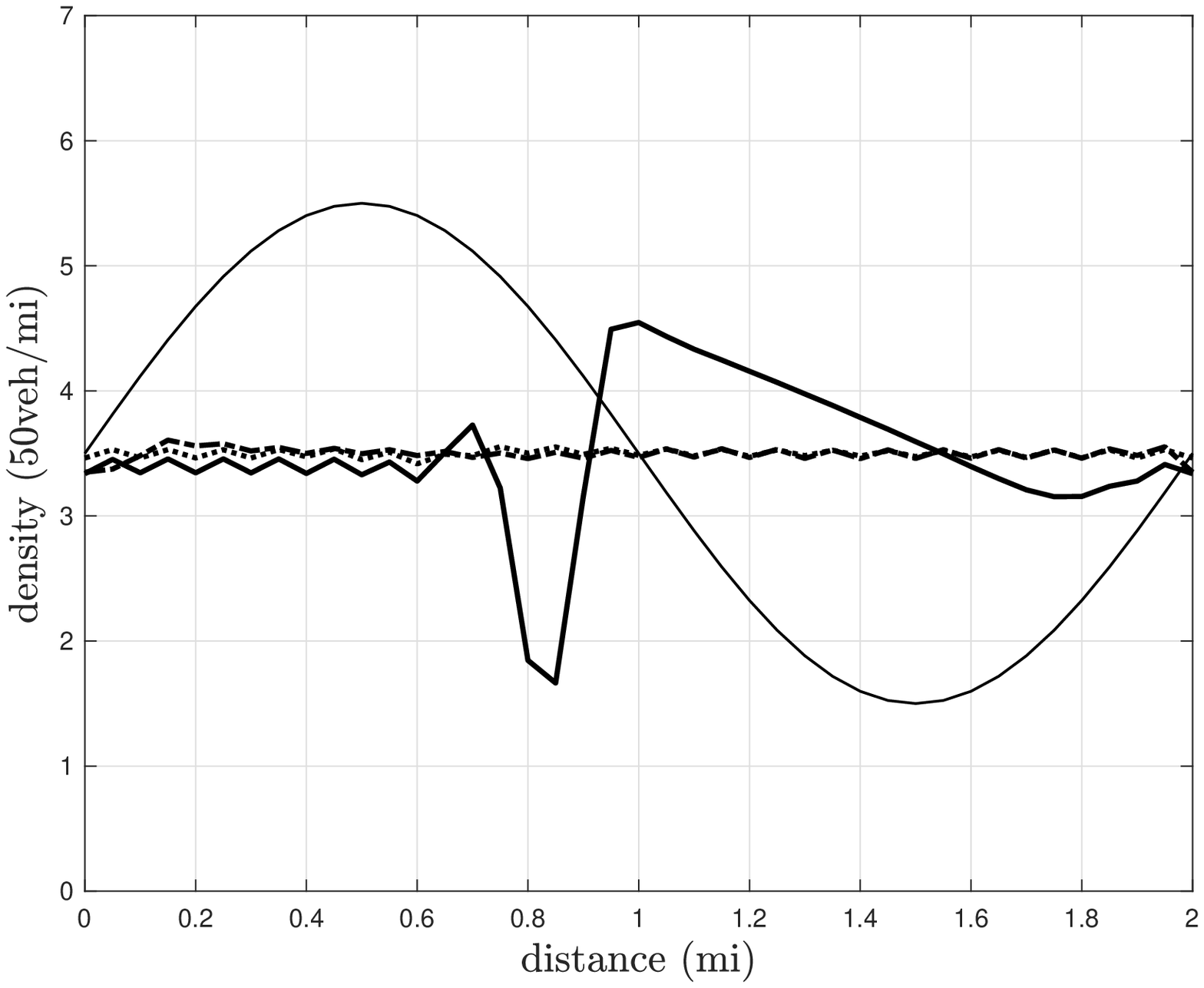}    
\caption{Total vehicle density at $t = 0$ (solid, light), $t = T/4$ (solid), $t = T/2$ (dot-dashed), and $t = T$ (dotted) corresponding to optimization of $\dot H^{-1}$ cost with $N=2$ and an upper limit on AV flow} 
\label{fig:rho_hist_H1_N2_lim}
\vspace{-0.2in}
\end{center}
\end{figure}

\begin{figure}[t]
\begin{center}
\includegraphics[width=0.45\textwidth]{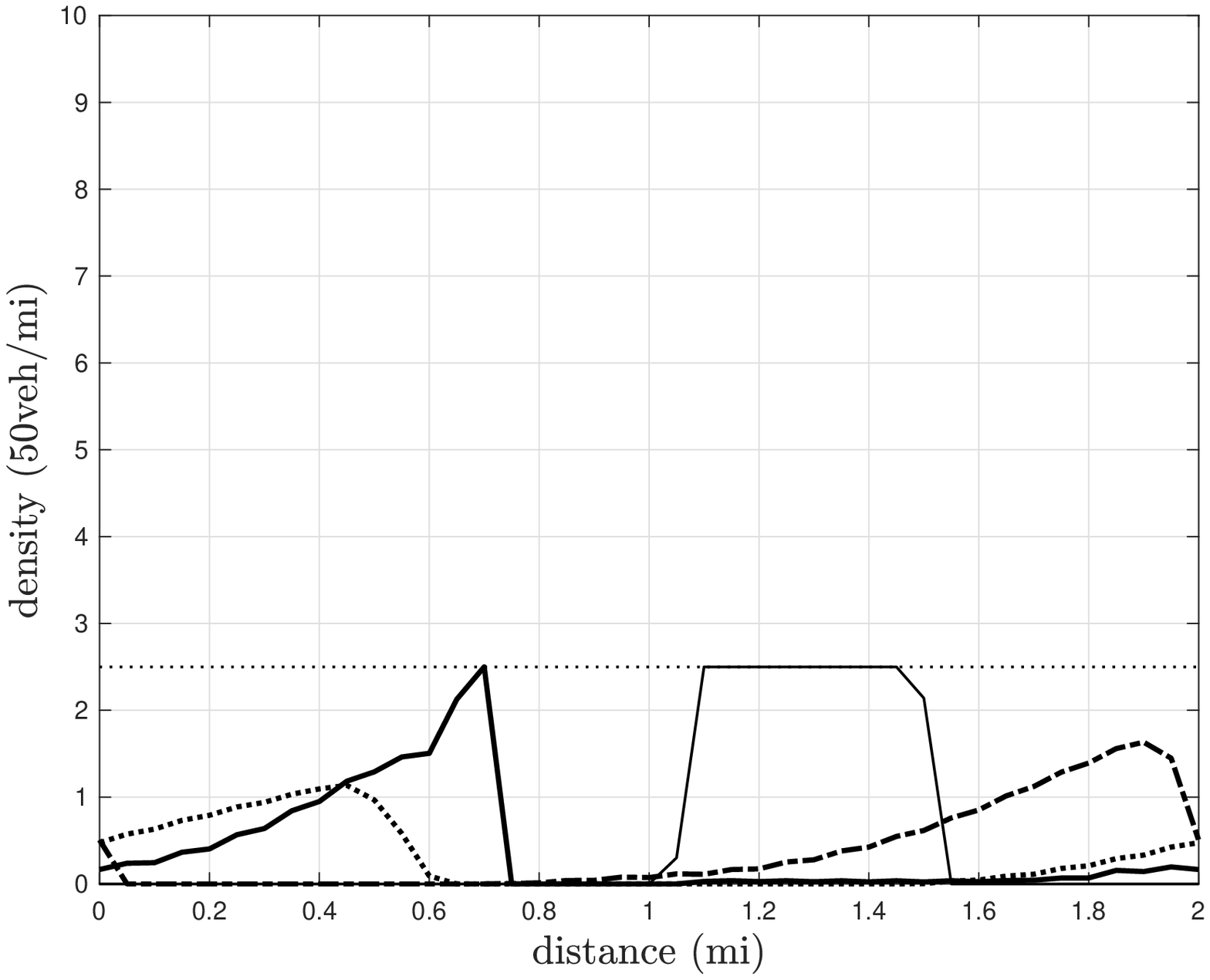}    
\caption{AV flow at $t = 0$ (solid, light), $t = T/4$ (solid), $t = T/2$ (dot-dashed), and $t = T$ (dotted) corresponding to optimization of $\dot H^{-1}$ cost with $N=2$ and an upper limit at $2.5$ (dotted, light)} 
\label{fig:H1_flow_lim}
\vspace{-0.2in}
\end{center}
\end{figure}

\subsection{Optimization of $L^2$ cost}

In the first simulation, we determine the control input according to the solution of the optimization problem \eqref{equ:L2_cost}-\eqref{equ:opt_constr}, which penalizes density according to the $L^2$ cost. We set the incentivization parameter $\beta = 0.2$ so that $0.2$ of vehicles can be AVs and we set $N=1$ so that we only solve the optimization once before implementing the control.

The time history is shown in Fig.~\ref{fig:rho_hist_L2_N1} and a qualitative evaluation suggests that we have not achieved a substantial improvement over the baseline. In particular, the density distribution at the final time is very close to the density distribution in the uncontrolled case. Shortening the control horizon by setting $N > 1$ achieves an even worse result; this is shown in Fig.~\ref{fig:rho_hist_L2_N2}, for which we set $N=2$. Investigating this further, we ran simulations for $N = 1,2,4,8,16$ for which we plot in Fig.~\ref{fig:L2_cost} the normalized time history of the density component of the cost \eqref{equ:L2_cost},
\begin{equation*}
\frac{1}{3.5^2|\Omega|}\int_{\Omega} (\rho_1(x,t)+\rho_2(x,t))^2 \text{d}x.
\end{equation*}
The plot shows that the we are able to initially, greatly reduce the density component of the cost in all cases but that this does not imply that we evenly distribute traffic in the steady state. In fact, for $N > 1$, we have worsened performance. For this reason, we consider the minimization of $\dot H^{-1}$.


\subsection{Optimization of $\dot H^{-1}$ cost}

The reason we believe that the $L^2$ cost optimization was not able to achieve good steady-state convergence is that being near to the optimal in terms of $L^2$ does not mean that the distribution is near uniform; due to the equivalence of norms, we expected that another choice of $L^p$ cost would have similar performance and this has been confirmed by simulations that are not reported here. An alternative is to consider the family of multiscale norms, typically used in fluid mixing. 

In the first simulation, we determine the control input according to the solution of the optimization problem \eqref{equ:H1_cost}, \eqref{equ:opt_constr}, which penalizes density according to the $\dot H^{-1}$ cost. We keep the incentivization parameter at $\beta = 0.2$ and set $N=1$.

In the results, not shown here, the density reaches close-to-uniform convergence at time $T/2$ but, by the end of the simulation, it drifts away from uniform. We suspect that the drift away from uniform is due to the fact that we optimize over a finite time-horizon, and so do not sufficiently penalize steady state. This is supported by the results of Figs.~\ref{fig:rho_hist_H1_N2}-\ref{fig:H1_cost}, for which we have set $N=2$ and obtained sustained convergence to uniform. Fig.~\ref{fig:rho_hist_H1_N2} shows sustained convergence to the uniform and Fig.~\ref{fig:H1_cost} shows that the drift away from uniform in the case of $N=1$ occurs toward the end of the simulation.

\subsection{Discussion}

Having obtained a satisfactory numerical results, we now discuss the practicality and implementability of this approach. In Fig.~
\ref{fig:rho2_hist_H1_N2}, we show the 
AV density, 
and in Fig.~\ref{fig:H1_flow}, we show the AV flow. The AV density shows that the optimal approach is to choose $0.2$ of the slower vehicles and switch these into autonomous mode, recomputing the distribution again halfway through the time-period.

From the plot of AV flow, we see that the required flow to obtain the desired result is higher than the maximum achievable in a single lane; to achieve it, it would require AVs to drive at the maximum rate in four separate lanes. This suggests the use of at least one dedicated lane for AVs in order to implement the scheme, along with a limit on the flow $m_2$. We perform an additional simulation, setting the limit at the maximum rate according to the LWR model \eqref{equ:LWR},
\begin{equation}
\frac{u_0\rho^*}{4} = 2.5,
\end{equation}
which corresponds to the use of at most one additional dedicated AV lane. 
The results are plotted in Figs.~\ref{fig:rho_hist_H1_N2_lim}-\ref{fig:H1_flow_lim}. Fig.~\ref{fig:H1_flow_lim} shows a plot of the AV flow at different times, exhibiting adherence to the prescribed limit. In Fig.~\ref{fig:rho_hist_H1_N2_lim}, the plot of total vehicle flow is not much qualitatively different than that of the unconstrained case, plotted in Fig.~\ref{fig:rho_hist_H1_N2}, validating our results. Note that we did not consider lane-change interactions between vehicles and assumed a bulk fluid model.

Before concluding, we note that the general scheme presented in this work is broadly applicable to most traffic scenarios and not limited to the single, circular road example on which we performed the numerical simulation. Nevertheless, this requires more future work as different scenarios will exhibit different topologies $\Omega$ and likely require a more specialized approach towards optimization.

\section{Conclusion}
\label{sec:conclude}

In this work, we presented an optimization-based control scheme for stabilizing traffic to a uniform state by incentivizing and controlling a group of autonomous vehicles. We considered two approaches, the first optimizing an $L^2$ cost and the second optimizing an $\dot H^{-1}$ cost.

We performed numerical simulations to determine the efficacy of the approach. The results showed that, by optimizing the $\dot H^{-1}$ cost, we are able to stabilize traffic to the uniform in finite time while, by optimizing the $L^2$ cost, we are not. We discussed the practicality of this approach and recommend the addition of a dedicated lane for implementability.

\addtolength{\textheight}{-12cm}   




\section*{Acknowledgment}

The authors acknowledge Dr.~Saleh Nabi of Mitsubishi Electric Research Laboratories for technical discussions.


\bibliographystyle{IEEEtran}
\bibliography{twopop}

\end{document}